\documentclass[12pt,oneside]{amsart}
\usepackage{calc,geometry, verbatim, graphicx, amssymb, color, mathtools,adjustbox}

\usepackage[pdftex,bookmarks,colorlinks,breaklinks]{hyperref} 
\hypersetup{linkcolor=blue,citecolor=red,filecolor=dullmagenta,urlcolor=darkblue} 
\newtheorem{theorem}{Theorem}[section]
\newtheorem{lemma}[theorem]{Lemma}

\newcommand{\Mod}[1]{\ (\mathrm{mod}\ #1)}

\newcommand{\cc}{\mathcal}
\newcommand\roundup[1]{\left\lceil#1\right\rceil}
\newcommand\rounddown[1]{\left\lfloor#1\right\rfloor}

\makeatletter
      \def\@setcopyright{}
      \def\serieslogo@{}
      \makeatother      
\begin{document}
   \author{Amin  Bahmanian}
   \address{Department of Mathematics,
  Illinois State University, Normal, IL USA 61790-4520}
\title[Symmetric Layer-Rainbow Colorations of Cubes]{Symmetric Layer-Rainbow Colorations of Cubes}

   \begin{abstract}  
Can we  color the $n^3$ cells of an $n\times n\times n$ cube $L$ with $n^2$ colors in such a way that each layer parallel to each face contains each color exactly once and that the coloring is symmetric so that $L_{ij\ell}=L_{j\ell i}=L_{\ell ij}$ for distinct $i,j,\ell \in  \{1,\dots,n\}$, and $L_{iij}=L_{jj i}, L_{iji}=L_{jij}, L_{ij j}=L_{jii}$ for  $i,j\in \{1,\dots,n\}$?
 Using transportation networks, we show that such a coloring is possible if and only if $n\equiv 0,2 \Mod 3$ (with two exceptions, $n=1$ and $n\neq 3$). Motivated by the designs of experiments, the study of these  objects (without symmetry) was initiated by Kishen and Fisher  in  the 1940's. These objects are also closely related to orthogonal arrays whose existence has been  extensively investigated, and  they are  natural three-dimensional analogues of symmetric latin squares. 

   \end{abstract}
   \subjclass[2010]{05B15, 05C70, 05C65, 05C15}
   \keywords{Latin cube, symmetry, factorization, circulation, transportation network, parallel class, orthogonal array}
   \date{\today}
   \maketitle

\section{Introduction}
A {\it $(2,1)$-latin cube} of {\it order} $n$ is an $n\times n\times n$ array on $n^2$ {\it symbols}, $\{1,2,\dots,n^2\}$, such that each layer parallel to each face contains each symbol exactly once. For example, stacking 
\begin{tabular}{ c c }
 1 & 2  \\ 
 3 & 4     
\end{tabular}
on top of  
\begin{tabular}{ c c }
 4 & 3  \\ 
 2 & 1     
\end{tabular}
forms a $(2,1)$-latin cube of order 2. The study of latin cubes was initiated in the 1940's by preeminent statisticians including  Fisher \cite{MR13113},   Kishen \cite {Kishen42,MR34743}, Rao \cite{MR19291}, and Brownlee and Loraine \cite{MR27485}. These objects  are also closely related to orthogonal arrays \cite{MR1693498, MR3495977} which themselves  have  connections to combinatorial designs, and have applications in  coding theory and cryptography \cite{MR2246267}. 
    
For $n,d,s,f\in \mathbb{N}$ with  $d\geq s+f$,  an {\it $(s,f)$-latin hypercube} of {\it order}  $n$ and {\it dimension} $d$ is an $n\times \dots \times n$ ($d$ times) array on $n^s$ symbols such that each  symbol occur exactly $n^{d-s-f}$ times in every subarray that  is obtained by fixing $f$ coordinates.  The existence of  $(s,f)$-hypercubes for given parameters is an open problem, and  there is no $(s,f)$-latin hypercube of order $n$ and dimension $d$ for $s\geq 2$ and $f>(n-1)^{s-1}$ (See \cite[Lemma 2.2]{MR2860603}.  
For further results on latin cubes, see \cite{TranCamBaiPraSch, MR2902643, MR3897549, MR3600882, MR2399374} and references therein.

Latin cubes with symmetric properties have been very little understood. For example, totally symmetric $(1,1)$-latin cubes of small order were enumerated by Bailey, Preece, and Zemroch \cite{MR511521}. Although there is no consensus on which one of $(1,1)$-latin, $(1,2)$-latin, and $(2,1)$-latin cubes should be called latin cubes, for the sake of convenience, in this paper, $(2,1)$-latin cubes will be called latin cubes. We introduce symmetric latin cubes and we settle their existence. A latin cube $L$ of order $n$ is {\it symmetric} if it meets the following conditions.
\begin{enumerate}
    \item $L_{ij\ell}=L_{j\ell i}=L_{\ell ij}$ for distinct $i,j,\ell \in  \{1,\dots,n\}$,
    \item and $L_{iij}=L_{jj i}, L_{iji}=L_{jij}, L_{ij j}=L_{jii}$ for  $i,j\in \{1,\dots,n\}$.
\end{enumerate}  
These objects are natural three-dimensional analogues of symmetric latin squares.  For $i\in \{1,\dots,n\}$, let $L_{i**}=\{L_{ij\ell}\ |\  1\leq j,\ell\leq n\}$, $L_{*i*}=\{L_{ji\ell}\ |\  1\leq j,\ell\leq n\}$, and $L_{**i}=\{L_{j\ell i}\ |\  1\leq j,\ell\leq n\}$. If $L$ is symmetric, and one of the three layers in $\{L_{i**},L_{*i*},L_{**i} \}$ is given, the other two layers can be uniquely determined. 

The main result of this paper is the following.
\begin{theorem} \label{symlcubeexithm}
There exists a symmetric latin cube of order $n$ if and only if  $n\equiv 0,2 \Mod 3$  with two exceptions: $n=1$ and $n\neq 3$. 
\end{theorem}

We show that the existence of a symmetric latin cube of order $n$ is equivalent to the existence of a partition  into parallel classes of  some non-uniform set system on $n$ points. Using transportation networks, we show that such a partition exists if and only if a certain system of equations has an integral solution. To complete the proof, we establish when the desired system has an integral solution.

\section{Symmetric Latin Cubes and Factorizations}

For a set $X$ and $a\in \mathbb N$, $a X$ is the multiset consisting of $a$ copies of each element of $X$, and $\binom{X}{a}$ is the collection of all $a$-subsets of $X$.   A partition $\cc P$ of a set $X$ is a  {\it parallel class}  of $\binom{X}{1}\cup 3\binom{X}{2} \cup 2\binom{X}{3}$ if $\cc P\subseteq\binom{X}{1}\cup \binom{X}{2} \cup \binom{X}{3}$.  Observe that if $\binom{X}{1}\cup 3\binom{X}{2} \cup 2\binom{X}{3}$ is  partitioned into parallel classes, then the number of parallel classes is
$$
1+ 3(|X|-1)+2\binom{|X|-1}{2}=|X|^2.
$$

\begin{lemma}
A symmetric latin cube of order $|X|$ exists if and only if $\binom{X}{1}\cup 3\binom{X}{2} \cup 2\binom{X}{3}$ can be partitioned into parallel classes.
\end{lemma}
\begin{proof}
Let $X=\{x_1,\dots, x_n\}$, and $\cc G=\binom{X}{1}\cup 3\binom{X}{2} \cup 2\binom{X}{3}$.

Suppose that there exists a symmetric latin cube $L$ of order $n$. We color each element of  $\cc G$ using one of the symbols in $\{1,\dots,n^2\}$ of $L$ in the following way. For $i\in \{1,\dots, n\}$, we color $\{x_i\}$ with $c$ if $L_{iii}=c$.   This lead to a coloring of $\binom{X}{1}$. Now, let $i,j\in \{1,\dots, n\}$ with $i\neq j$. Suppose $L_{iij}=L_{jji}=c, L_{iji}=L_{jij}=c', L_{jii}=L_{ijj}=c''$. Since $L$ is latin, $c,c',c''$ are pairwise distinct. We color a copy of $\{x_i,x_j\}$ with $c$, a copy with $c'$ and a copy with $c''$.  Finally, let $i,j,\ell$ be three distinct elements of $\{1,\dots, n\}$.  Suppose $L_{ij\ell}=L_{j \ell i}=L_{\ell i j}=c, L_{ji \ell}=L_{i \ell j}=L_{\ell j i}=c'$. Again, since $L$ is latin, $c$ and $c'$ are distinct. We color a copy of $\{x_i,x_j, x_\ell\}$ with $c$, and a copy with $c'$.  Let us verify that each color class is a parallel class of $\cc G$. Fix $x_i\in X$. The sets in $\cc G$ that contain $x_i$ correspond to the cells in layer $L_{i**}$. Since $L$ is latin, the $n^2$ symbols in $L_{i**}$ are distinct, and consequently, the color of the $n^2$ sets containing $x_i$  are different. 

Conversely, suppose that $\cc G$ can be partitioned into parallel classes. First, fill $L_{iii}$ with the color of  $\{x_i\}$ for  $i\in \{1,\dots, n\}.$ Suppose the color of three three copies of $\{x_i,x_j\}$ are $c,c'$ and $c''$. Since these three colors are  distinct, we fill  the cells $L_{iij}$ and $L_{jji}$ with $c$, $L_{iji}$ and $L_{jij}$ with $c'$, and $L_{jii}$ and $L_{ijj}$ with $c''$. 
Now, suppose the color of the two copies of $\{x_i,x_j, x_{\ell}\}$ are $c$ and $c'$ (which are distinct). We fill the cells  $L_{ij\ell}$, $L_{j\ell i}$ and $L_{\ell ij}$ with $c$, and $L_{i\ell j}$, $L_{\ell j i}$, and $L_{ji \ell}$ with $c'$. It is clear that $L$ is symmetric. Because of the one-to-one correspondence between the sets in $\cc G$ that contain $x_i$ and the cells in layer $L_{i**}$, $L$ is latin.
\end{proof}

\subsection{An Example}
\begin{figure}[p] 
  \begin{adjustbox}{addcode={\begin{minipage}{\width}}{\caption{
      A symmetric latin cube of order five together with all its layers
      }\end{minipage}},rotate=90,center}  
 \includegraphics[width=10in,height=6in,keepaspectratio]{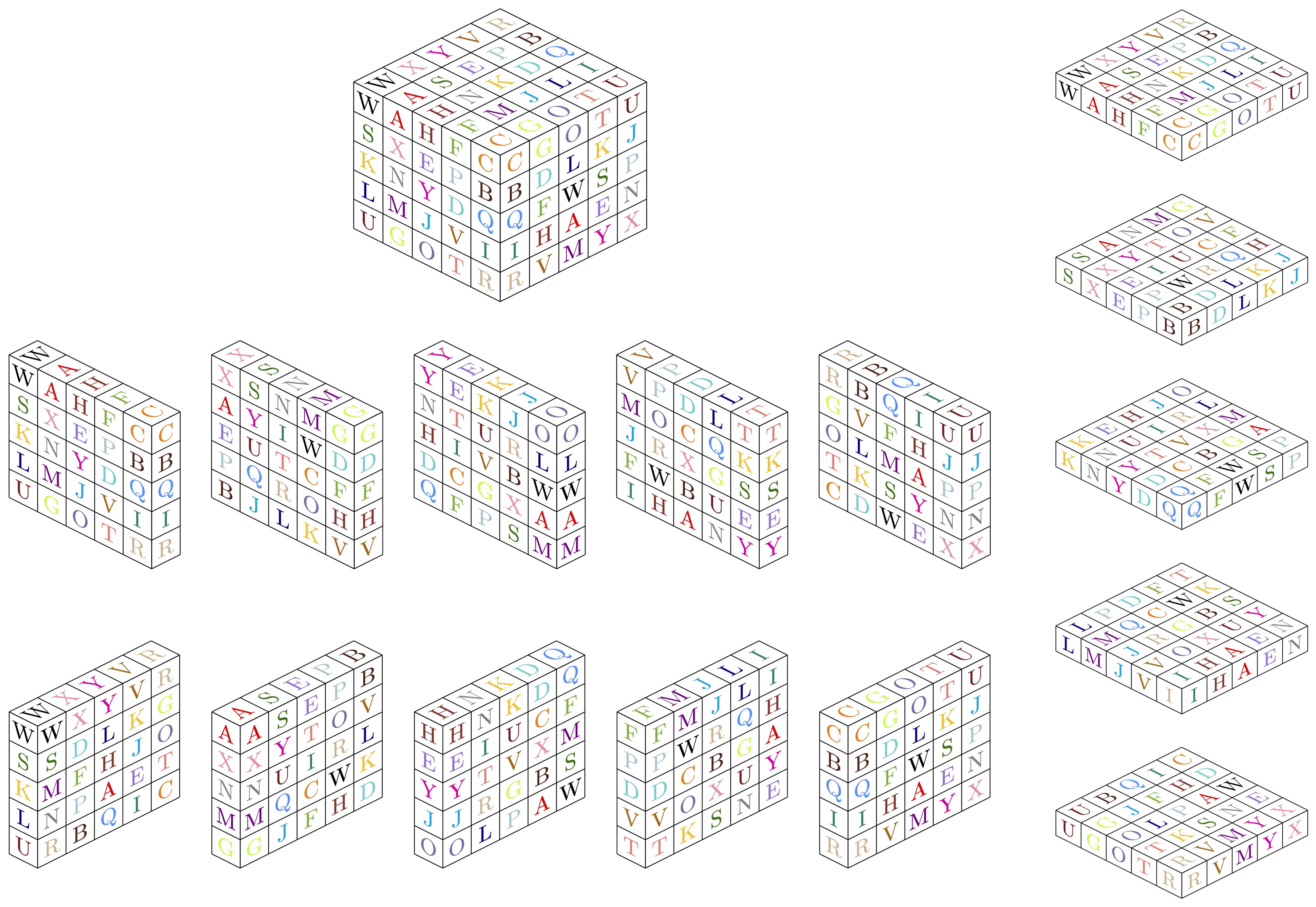} \end{adjustbox}
 \label{fig5symlatincube} 
\end{figure}

Let $X=\{1,2,3,4,5\}$. Consider the following partition of $\binom{X}{1}\cup 3\binom{X}{2} \cup 2\binom{X}{3}$ into parallel classes (We abbreviate sets such as $\{x,y,z\}$ to $xyz$). 
\begin{align*}
\big\{&\{12, 345\},\{ 34, 125\}, \{15, 234\}, \{ 25, 134\},\{ 45, 123\}, \{ 14, 235\}, \{34, 125\}, \{13, 245\}, \\
&\{23, 145\}, \{25, 134\}, \{13, 245\}, \{14, 235\}, \{35, 124\}, \{45, 123\}, \{24, 135\}, \{35, 124\}, \\
&\{24, 135\}, \{15, 234\}, \{12, 345\}, \{23, 145\}, \{4, 15, 23\}, \{3, 14, 25\}, \{1, 24, 35\}, \\
&\{5, 12, 34\}, \{2, 13, 45\}\big\}.
\end{align*}
Assigning symbols $A, B,\dots, Y$ to the 26 parallel classes, respectively, we obtain the following  symmetric latin cube of order 5 (See Figure \ref{fig5symlatincube}). 

\begin{center}
\begin{tabular}{ |c |c| c| c| c|}
\hline  W & A & H & F & C \\ \hline
S & X & E & P & B \\ \hline
K & N & Y & D & Q \\ \hline
L & M & J & V & I \\ \hline
U & G & O & T & R \\ \hline
\end{tabular}\quad 
\begin{tabular}{ |c |c| c| c| c|}\hline
 X & S & N & M & G \\ \hline
A & Y & I & W & D \\ \hline
E & U & T & C & F \\ \hline
P & Q & R & O & H \\ \hline
B & J & L & K & V \\ \hline
\end{tabular}\quad 
\begin{tabular}{ |c |c| c| c| c|}\hline
 Y & E & K & J & O \\ \hline
N & T & U & R & L \\ \hline
H & I & V & B & W \\ \hline
D & C & G & X & A \\ \hline
Q & F & P & S & M \\ \hline
\end{tabular}
\end{center}
\begin{center}
\begin{tabular}{ |c |c| c| c| c|}\hline
 V & P & D & L & T \\ \hline
M & O & C & Q & K \\ \hline
J & R & X & G & S \\ \hline
F & W & B & U & E \\ \hline
I & H & A & N & Y \\ \hline
\end{tabular}\quad 
\begin{tabular}{ |c |c| c| c| c|}\hline
 R & B & Q & I & U \\ \hline
G & V & F & H & J \\ \hline
O & L & M & A & P \\ \hline
T & K & S & Y & N \\ \hline
C & D & W & E & X \\ \hline
\end{tabular}
\end{center}

\section{A System  of Equations}
\begin{lemma} \label{syslem}
For $n\in \mathbb N$,  the following system has a nonnegative integral solution if and only if  $n\equiv 0,2 \Mod 3$  with two exceptions: $n=1$ and $n\neq 3$. 
\begin{align} \label{system1}
    \begin{cases}
       a_i+2b_i+3c_i=n &\mbox { for } i\in \{1,\dots, n^2\},\\ 
        \displaystyle\sum_{i=1}^{n^2} a_i=n,\\
        \displaystyle\sum_{i=1}^{n^2} b_i=3\binom{n}{2},\\
        \displaystyle\sum_{i=1}^{n^2} c_i=2 \binom{n}{3}.
    \end{cases}
\end{align}
\end{lemma}
\begin{proof}
If $n=1$, then $a_1=1,b_1=c_1=0$ is the only solution to \eqref{system1}. If on the contrary  \eqref{system1} has a solution for $n=3$, then since $a_i+2b_i+3c_i=3$, we  have $(a_i,b_i,c_i)\in \{(3,0,0), (1,1,0), (0,0,1)\}$ for $i\in \{1,\dots, 9\}$. But $\sum_{i=1}^9 b_i=9$,  so we  have  $b_1=\dots=b_9=1$, and so $a_1=\dots=a_9=0$ which contradicts $\sum_{i=1}^9 a_i=3$. 

Now, let $n>1, n\equiv 1\Mod 3$. Suppose on the contrary that \eqref{system1} has a solution. Since $a_i+2b_i+3c_i=n$, we have $a_i+2b_i\equiv 1 \Mod 3$, and consequently, either $a_i\geq 1$ or $b_i\geq 2$ for $i\in \{1,\dots, n^2\}$. Let $A$ be the set of indices $i\in \{1,\dots, n^2\}$ with $a_i\geq 1$, and let $B$ be the set of indices $i\in \{1,\dots, n^2\}$ with $a_i=0$ and $b_i\geq 2$. Note that $|A\cup B|=n^2$, $|A|\leq n$ and $|B|\leq \frac{3}{2}\binom{n}{2}$. Since $A\cap B=\emptyset$, we have
$n^2\leq n+\frac{3}{2}\binom{n}{2}$, or equivalently, $n\leq 1+\frac{3}{4}(n-1)$. Therefore, $n\leq 1$, which is a contradiction.

To solve \eqref{system1} in the remaining cases, we find a $3\times n^2$ nonnegative integral matrix $M_n$ whose row-sums are $n, 3\binom{n}{2}$, and $2\binom{n}{3}$ such that 
\begin{align*}
   (1\  2 \  3) M_n=(n\ \dots\  n)_{1\times n^2}.
 \end{align*}
 In each case, it is a routine exercise to check that the desired properties are satisfied. Here, $\mathbf{0}$ is a block of zeros and  $J_t:=(1\  \dots\  1)_{1\times t}$ for $t\in \mathbb N$. 
\begin{enumerate} 
    \item [(a)] $n\equiv 0\Mod 6$: Let 
$$
M_n=\left(
\begin{array}{ccc}
{n} & \mathbf{0}   & \mathbf{0}\\
0 &\dfrac{n}{2} J_{3(n-1)} &  \mathbf{0}\\
0 & \mathbf{0} &  \dfrac{n}{3} J_{(n-1)(n-2)}
\end{array}
\right).
$$   
\item [(b)] $n\equiv 2\Mod 6$: Let 
$$
M_n=\left(
\begin{array}{ccc}
{n} & \mathbf{0}   & \mathbf{0}\\
0 &\dfrac{n}{2} J_{n-1} &  J_{n(n-1)}\\
0 & \mathbf{0} &  \dfrac{n-2}{3} J_{n(n-1)}
\end{array}
\right).
$$   

\item [(c)] $n\equiv 3\Mod 6, n\neq 3$:  If $n\geq 27$,  let
 $$ M_n=
\left(
\begin{array}{cccc}
{n} & \mathbf{0}   & 0 & \mathbf{0}\\
0 &\dfrac{n-3}{2} J_{3(n+2)} & 9 &  \mathbf{0}\\
0 & J_{3(n+2)}  & \dfrac{n-18}{3} & \dfrac{n}{3} J_{n^2-3n-8}
\end{array}
\right).
$$    
For $n\in \{9,15,21\}$,  let 
\begin{align*}
M_9=\left(
\begin{array}{ccc}
{9} &  \mathbf{0}  & \mathbf{0}\\
0 &3 J_{36} &  \mathbf{0}\\
0 & J_{36} &  3J_{44}
\end{array}
\right), 
M_{15}=\left(
\begin{array}{cccc}
{15} & \mathbf{0}   & 0& \mathbf{0}\\
0 &6 J_{52} &  3&\mathbf{0}\\
0 & J_{52} & 1 & 5J_{171}
\end{array}
\right), 
M_{21}=\left(
\begin{array}{ccc}
{21} &  \mathbf{0}  & \mathbf{0}\\
0 &9 J_{70} &  \mathbf{0}\\
0 & J_{70} &  7J_{370}
\end{array}
\right).
\end{align*}
\item [(d)] $n\equiv 5\Mod 6$: If $n\geq 17$,  let
 $$
M_n=\left(
\begin{array}{cccc}
{n} & 0   & \mathbf{0} & \mathbf{0}\\
0 &7 &  J_{n+2} &  \dfrac{n-3}{2} J_{n^2-n-4}\\
0 & \dfrac{n-14}{3}  & \dfrac{n-2}{3}J_{n+2} & J_{n^2-n-4}
\end{array}
\right).
$$    
For $n\in \{5,11\}$,  let 
\begin{align*}
M_5=    \left(
\begin{array}{cc}
  \mathbf{0}  & J_5\\
 J_{20} &  2J_{5}\\
 J_{20} &  \mathbf{0}
\end{array}
\right),\quad 
M_{11}=
\left(
\begin{array}{ccc}
{11} &  \mathbf{0}  & \mathbf{0}\\
0 & J_{105} &  4J_{15}\\
1 & 3J_{105} &  J_{15}
\end{array}
\right).
\end{align*}
\end{enumerate}
\end{proof}   
\section{Circulations}
A {\it circulation} on a digraph $D$ is a  real-valued function $f$ on the edge set of $D$ which satisfies the  conservation condition at every vertex, i.e. at each vertex, the sum of the values of $f$ on incoming edges equals the sum of the values of $f$ on outgoing edges. By  $(v,w)$ we mean a directed edge  from $v$ to $w$, and we abbreviate $f(\{v,w\})$ to  $f(v,w)$.

For a multiset $A$ and $u\in A$, let $\mu_A (u)$ denote the multiplicity of $u$ in $A$, and $A=\{u_1^{r_1},\dots,u_k^{r_k}\}$ is a multi-set in which $\mu_A (u_i)=r_i$ for $i\in \{1,\dots,k\}$. For multisets $A_1,\dots, A_n$, we define $U=\bigcup_{i=1}^n A_i$ so that $\mu_{U}(u)=\sum_{i=1}^n \mu_{A_i}(u)$ for $u\in U$; for example $\{u^2,v,w^2\}\cup \{u,w^2\}=\{u^3,v,w^4\}$. A {\it split} of $\{1,\dots,n\}$ is a multi-set $\mathcal F$ such that each $A\in \mathcal F$ is of the form $\{1^r\}\cup S$ for some $0\leq r\leq 3$ and for a set $S\subseteq \{2,\dots, n\}$ with $0\leq |S|\leq 3-r$. For example, $\mathcal F=\{\{1,2,3\}, \{1\}, \{1^2\}, \{1,4\}, \{1^2,5\}\}$ is a split of $\{1,2,3,4,5\}$.

Let $\ell \in \{1,\dots, n\}$. For $i\in \{1,\dots, n^2\}$, $j\in \{1,\dots, \ell\}$  and a multiset $\{\cc F_1,\dots, \cc F_{n^2}\}$ of splits of $\{1,\dots, \ell\}$, let $\mu_i^j=\sum_{A\in \cc F_i}\mu_A(j)$. For $0\leq r\leq 3$ and a set  $S\subseteq \{2,\dots, \ell\}$ with $0\leq |S|\leq 3-r$, let $\mu_S^r=\sum_{i=1}^{n^2} \mu_{\cc F_i}\big( S\cup \{1^r\}\big)$, and let
\begin{align*}
    a_S^r=\begin{cases}
          0 & \mbox { if } r+|S|=0,\\
          1 & \mbox { if } r+|S|=1,\\
          3 & \mbox { if } r+|S|=2,\\
          2 & \mbox { if } r+|S|=3.
    \end{cases}
\end{align*}
The case $\ell=n$ of the following results settles Theorem \ref{symlcubeexithm} (See Theorem \ref{finalsymlatcubthm}). 
\begin{theorem}\label{seemstmainres}
Let $n\in \mathbb N$, $n\neq 3$, and  $n\equiv 0,2 \Mod 3$. For any  $\ell \in \{1,\dots, n\}$, there exists a multiset $\{\mathcal F_1,\dots,\mathcal F_{n^2}\}$ of splits of $\{1,\dots,\ell\}$ such that the following conditions hold.
\begin{align*}
  \begin{cases}
        \mu_i^1=n-\ell+1 &\mbox { for } i\in \{1,\dots, n^2\},\\
        \mu_i^j=1 & \mbox { for } i\in \{1,\dots, n^2\}, j\in \{2,\dots,\ell\},
        \\ 
    \mu_S^r=a_S^r\dbinom{n-\ell+1}{r} &\mbox { for } 0\leq r\leq 3, \mbox { and for a set } S\subseteq \{2,\dots,\ell\} \mbox{ with } 0\leq |S|\leq 3-r.
\end{cases}\end{align*}
\end{theorem}
\begin{proof}
The proof is by induction on $\ell$. By Lemma \ref{syslem}, there exit non-negative integers $a_i,b_i,c_i$ for $i\in \{1,\dots, n^2\}$ such that \eqref{system1} holds. 
For $\ell=1$, we define the following splits.
$$\cc F_i=\Big\{\{1\}^{a_i}, \{1^2\}^{b_i}, \{1^3\}^{c_i}\Big\}\quad \mbox{ for } i\in \{1,\dots, n^2\}.$$
By the following, $\{\mathcal F_1,\dots,\mathcal F_{n^2}\}$ satisfies the desired conditions. 
\begin{align*}
\begin{cases}
  \mu_i^1=a_i+2b_i+3c_i=n &\mbox { for } i\in \{1,\dots, n^2\},\\
    \mu_\emptyset^r=a_\emptyset^r\dbinom{n}{r} &\mbox { for } 0\leq r\leq 3.
\end{cases}
\end{align*}

Assume that for some value $\ell<n$, a multiset $\{\mathcal F_1,\dots,\mathcal F_{n^2}\}$ of splits of $\{1,2\dots,\ell\}$ with the required properties exist. We form a digraph $D$ with vertex multiset 
\begin{align*}
    V&=\Big\{\sigma, \tau\Big\}\\
    &\  \bigcup \Big\{v_1,\dots,v_{n^2}\Big\} \\
    &\  \bigcup \Big\{v_{i}^A \ \Big |\  i\in \{1,\dots, n^2\}, \ A\in  \mathcal F_i, \mu_A(1)>0 \Big\}\\ 
    &\  \bigcup \Big\{v_S^r \ \Big |\ 1\leq r\leq 3, S\subseteq \{2,\dots, \ell\}, 0\leq |S| \leq 3-r \Big\},
\end{align*}    
(since each $\cc F_i$ is a multiset, there may be multiple occurrences of $v_{i}^A$ in $V$) and with a circulation $f$ as follows.  
\begin{itemize}
    \item For $i\in \{1,\dots,n^2\}$, there is a directed edge from $\sigma$ to $v_i$ such that $f(\sigma,v_i)=1$.
    \item For $i\in \{1,\dots,n^2\}$ and for $A\in \mathcal F_i$ with $\mu_A(1)>0$, there is a directed edge from $v_i$ to (each occurrence of) $v_i^A$ such that $f(v_i,v_i^A)=\mu_A(1)/(n-\ell+1)$.
    \item For $i\in \{1,\dots,n^2\}$ and for $A\in \mathcal F_i$ with $\mu_A(1)>0$, if $A=S\cup \{1^r\}$, there is a directed edge from $v_i^A$ to $v_S^r$ such that $f(v_i^A,v_S^r)=\mu_A(1)/(n-\ell+1)=r/(n-\ell+1)$.
    \item For $1\leq r \leq 3$, and $S\subseteq \{2,\dots, \ell\}$ with $0\leq |S| \leq 3-r$, there is a directed edge from $v_S^r$ to $\tau$ such that $f(v_S^r,\tau)=a_S^r \binom{n-\ell}{r-1}$.
    \item There is a directed edge from $\tau$ to $\sigma$ such that $f(\tau, \sigma)=n^2$.
\end{itemize}
Let us verify that $f$ is a circulation. It is clear that $f$ satisfies the  conservation condition at $\sigma$ and at each $v_i^A$. By the induction hypothesis, $\mu_i^1=n-\ell+1$ for $i\in \{1,\dots,n^2\}$, and so the  conservation condition holds at each $v_i$. At each $v_S^r$, the sum of the values of $f$ on incoming edges equals (Here, $A=S\cup\{1^r\}$)
$$
\sum_{i=1}^{n^2}\sum_{\substack{A\in \cc F_i, \\ \mu_A(1)=r}} \frac{r}{n-\ell+1} = \frac{r}{n-\ell+1} \mu_{A\backslash \{1^r\}}^r=\frac{r a_S^r}{n-\ell+1} \binom{n-\ell+1}{r}= a_S^r\binom{n-\ell}{r-1},
$$
and so  the  conservation condition holds at each $v_S^r$. Finally, at $\tau$, the sum of the values of $f$ on incoming edges equals
\begin{align*}
    &1+3(\ell-1)+2\binom{\ell-1}{2}+3(n-\ell)+2(n-\ell)(\ell-1)+(n-\ell)(n-\ell-1)=n^2.
\end{align*}

Since $f$ is a circulation, there exists an integer circulation $g$ on $D$ such that $g(e)\in \{\rounddown{f(e)}, \roundup{f(e)}\}$  for every edge $e$ (see \cite{MR2729968} or \cite[Theorem 7.4]{MR1871828}). Let us fix an $i\in \{1,\dots,n^2\}$. For each $A\in \mathcal F_i$ with $\mu_A(1)>0$, we have $g(v_i, v_i^A)\in \{0,1\}$, but since $g(\sigma, v_i)=1$,  there is exactly one $A$ in $\cc F_i$ such that $g(v_i, v_i^A)=1$.  
Now, we obtain a split $\bar {\cc F}_i$ of the set $\{1,\dots,\ell+1\}$ by letting $\bar {\cc F}_i$ be obtained from $\cc F_i$ by replacing one 1 in $A\in \cc F_i$ with $\ell+1$ if $g(v_i, v_i^A)=1$.

We show that the splits $\overline {\cc F}_1,\dots, \bar {\cc F}_{n^2}$ satisfy the required properties. We define $\bar \mu_i^j$ and $\bar \mu_S^r$ for $\bar {\cc F}_1,\dots, \bar {\cc F}_{n^2}$ similar to the way we defined them for  $\cc F_1,\dots,\cc F_{n^2}$. Obviously, $\bar \mu_i{^j}=\mu_i^j=1$ for $i\in \{1,\dots, n^2\}, j\in \{2,\dots,\ell\}$. 
Moreover, $\bar \mu^{\ell+1}_i=1$, 
and $\bar \mu^1_i=\mu_i^1-\bar \mu^{\ell+1}_i=n-\ell$ for $i\in \{1,\dots, n^2\}$. 
Now, let $0\leq r \leq 3$ and $S\subseteq \{2,\dots, \ell+1\}$ with $0\leq |S| \leq 3-r$. If $\ell+1\in S$, we have
$$\bar \mu_S^r=g(v_{S\backslash\{\ell+1\}}^{r+1},\tau)=f(v_{S\backslash\{\ell+1\}}^{r+1},\tau)=a_S^r \binom{n-\ell}{r},$$
and if $\ell+1\notin S$, we have
$$\bar \mu_S^r=\mu_S^r-g(v_S^r,\tau)=a_S^r\dbinom{n-\ell+1}{r}-a_S^r \binom{n-\ell}{r-1}=a_S^r\dbinom{n-\ell}{r}.$$
\end{proof}

\section{Proof of the Main Result}
\begin{theorem} \label{finalsymlatcubthm}
The collection $\binom{X}{1}\cup 3\binom{X}{2} \cup 2\binom{X}{3}$ can be partitioned into parallel classes  if and only if  $|X|\equiv 0,2 \Mod 3$  with two exceptions: $|X|=1$ and $|X|\neq 3$. 
\end{theorem}
\begin{proof}
First, suppose that $\binom{X}{1}\cup 3\binom{X}{2} \cup 2\binom{X}{3}$ can be partitioned into parallel classes. It is clear that \eqref{system1} has a nonnegative integral solution, or equivalently,  $|X|\equiv 0,2 \Mod 3$.

Conversely, suppose that $|X|\equiv 0,2 \Mod 3$. If we let $\ell=n:=|X|$ in Theorem \ref{seemstmainres},  we obtain a multiset $\{\mathcal F_1,\dots,\mathcal F_{n^2}\}$ of splits of $\{1,\dots,n\}$ satisfying the following conditions.
\begin{align}
       & \mu_i(j)=1  \mbox { for } i\in \{1,\dots, n^2\}, j\in \{1,\dots,n\}, \label{degcondln}\\
    & \mu_S^r=a_S^r\dbinom{1}{r} \mbox { for } 0\leq r\leq 3 \mbox { and a set } S\subseteq \{2,\dots,n\} \mbox{ with } 0\leq |S|\leq 3-r. \label{multcondln}
\end{align}
Condition \eqref{degcondln} implies that for $i\in \{1,\dots, n^2\}$, each split $\mathcal F_i$ in indeed a partition of $\{1,\dots,n\}$, so each $\cc F_i$ is a parallel class of $\binom{X}{1}\cup 3\binom{X}{2} \cup 2\binom{X}{3}$. Now, let $0\leq r\leq 3$, and $S\subseteq \{2,\dots,n\}$  with $0\leq |S|\leq 3-r$. By \eqref{multcondln}, $\mu_S^r=0$ for $r>1$, $\mu_S^0=a_S^0$ and $\mu_S^1=a_S^1$, and consequently, $\{\mathcal F_1,\dots,\mathcal F_{n^2}\}$ is a partition of $\binom{X}{1}\cup 3\binom{X}{2} \cup 2\binom{X}{3}$ into parallel classes.
\end{proof}

\bibliographystyle{plain}

\section{Appendix}
A symmetric latin cube of order 6:

\tabcolsep=0.055cm
\begin{tabular}{|l|l|l|l|l|l|}
\hline
\begin{tabular}{ c c c c c c}
 j & D & B & F & K & A \\ 
 C & M & U & h & P & W \\  
 L & V & H & d & T & Y  \\
 O & R & a & I & i & c \\
 N & Z & b & Q & J & X \\
 E & f & S & g & e & G 
\end{tabular} 
&
\begin{tabular}{ c c c c c c}
 M & C & V & R & Z & f \\ 
 D & j & A & J & B & L \\  
 U & I & O & X & g & Q  \\
 h & N & e & G & S & T \\
 P & H & c & Y & E & a \\
 W & L & i & b & d & K 
\end{tabular}
&  
\begin{tabular}{ c c c c c c}
 H & U & L & a & b & S \\ 
 V & O & I & e & f & i \\  
 B & A & j & M & C & D  \\
 d & X & E & K & f & Z \\
 T & g & F & W & G & h \\
 Y & Q & J & P & R & N 
\end{tabular}
& 
\begin{tabular}{ c c c c c c}
 I & h & d & O & Q & g \\ 
 R & G & X & N & Y & b \\  
 a & e & K & E & W & P  \\
 F & J & M & j & D & C \\
 i & S & f & A & L & U \\
 c & T & Z & B & V & H 
\end{tabular}
&  
\begin{tabular}{ c c c c c c}
 J & P & T & i & N & e \\ 
 Z & E & g & S & H & d \\  
 b & c & G & f & F & R  \\
 Q & Y & W & L & A & V \\
 K & B & C & D & j & O \\
 X & a & h & U & I & M 
\end{tabular}
&   
\begin{tabular}{ c c c c c c}
 G & W & Y & c & X & E \\ 
 f & K & Q & T & a & L \\  
 S & i & N & Z & h & J  \\
 g & b & P & H & U & B \\
 e & d & R & V & M & I \\
 A & F & D & C & O & j 
\end{tabular}
\\ \hline
\end{tabular}

\vspace{.3in}
A symmetric latin cube of order 8:
\begin{center}
\tabcolsep=0.055cm
\begin{tabular}{|l|l|l|l|}
\hline
\begin{tabular}{ c c c c c c c c }
 $\beta$ & Q & $\alpha$ & I & K & w & i & S \\ 
 G & 8 & h & M & n & X & t & C \\  
 o & g & L & W & y & H & q & V \\
 2 & s & x & 3 & p & P & v & Y \\
 4 & $\ell$ & d & e & E & r & u & c \\
 0 & A & 9 & R & b & 5 & O & F \\
 6 & N & B & J & Z & T & j & a \\
 1 & K & U & z & 7 & f & m & D
\end{tabular}
&
\begin{tabular}{ c c c c c c c c }
 8 & G & g & s & $\ell$ & A & N & k \\ 
 Q & $\beta$ & R & q & O & z & U & J \\  
 H & 1 & 6 & w & S & m & Y & K \\
 M & T & b & F & j & c & y & 9 \\
 n & a & f & i & H & B & W & v \\
 X & 3 & u & E & I & 7 & D & d \\
 t & 0 & P & V & o & p & 4 & e \\
 C & 2 & Z & L & x & $\alpha$ & r & 5
\end{tabular}
&
\begin{tabular}{ c c c c c c c c }
 L & h & o & x & d & 9 & B & U \\ 
 g & 6 & 1 & b & f & u & P & Z \\  
 $\alpha$ & R & $\beta$ & r & N & e & c & s \\
 W & w & 0 & 5 & D & $\ell$ & 7 & O \\
 y & S & 2 & A & M & C & F & T \\
 H & m & n & a & J & v & Q & i \\
 q & Y & 3 & G & z & 8 & k & I \\
 V & K & 4 & X & t & j & E & p
\end{tabular}
&
\begin{tabular}{ c c c c c c c c }
 3 & M & W & 2 & e & R & J & z \\ 
 s & F & w & T & i & E & V & L \\  
 x & b & 5 & 0 & A & a & G & X \\
 I & a & r & $\beta$ & g & Z & d & u \\
 p & j & D & 1 & m & U & $\alpha$ & 8 \\
 P & c & $\ell$ & 4 & k & t & S & o \\
 v & y & 7 & f & h & K & C & H \\
 Y & 9 & O & 6 & Q & N & n & B
\end{tabular}
\\ \hline
\begin{tabular}{ c c c c c c c c }
 E & n & y & p & 4 & b & Z & 7 \\ 
 $\ell$ & H & S & j & a & I & o & x \\  
 d & f & M & D & 2 & J & z & t \\
 e & i & A & m & 1 & k & h & Q \\
 K & O & N & g & $\beta$ & V & X & P \\
 r & B & C & U & 6 & Y & s & G \\
 u & W & F & $\alpha$ & 5 & L & 9 & w \\
 c & v & T & 8 & 0 & q & R & 3
\end{tabular}
&
\begin{tabular}{ c c c c c c c c }
 5 & X & H & P & r & 0 & T & f \\ 
 A & 7 & m & c & B & 3 & p & $\alpha$ \\  
 9 & u & v & $\ell$ & C & n & 8 & j \\
 R & E & a & t & U & 4 & K & N \\
 b & I & J & k & Y & 6 & L & q \\
 w & z & e & Z & V & $\beta$ & x & h \\
 O & D & Q & S & s & 1 & 2 & g \\
 F & d & i & o & G & y & M & W
\end{tabular}
&
\begin{tabular}{ c c c c c c c c }
 j & t & q & v & u & O & 6 & m \\ 
 N & 4 & Y & y & W & D & 0 & r \\  
 B & P & k & 7 & F & Q & 3 & E \\
 J & V & G & C & $\alpha$ & S & f & n \\
 Z & o & z & h & 9 & s & 5 & R \\
 T & p & 8 & K & L & 2 & 1 & M \\
 i & U & c & d & X & x & $\beta$ & b \\
 a & e & I & H & w & g & $\ell$ & A
\end{tabular}
&
\begin{tabular}{ c c c c c c c c }
 D & C & V & Y & c & F & a & 1 \\ 
 k & 5 & K & 9 & v & d & e & 2 \\  
 U & Z & p & O & T & i & I & 4 \\
 z & L & X & B & 8 & o & H & 6 \\
 7 & x & t & Q & 3 & G & w & 0 \\
 f & $\alpha$ & j & N & q & W & g & y \\
 m & r & E & n & R & M & A & $\ell$ \\
 S & J & s & u & P & h & b & $\beta$
\end{tabular}
\\ \hline
\end{tabular}
\end{center}

\end{document}